\def \dint {\displaystyle\int}
\def \dsum {\displaystyle\sum}

\documentclass[a4paper,12pt]{article}
\usepackage{amsfonts}
\usepackage{amsmath}
\newtheorem{theorem}{Theorem}

\begin{document}

\title{A New Generalization of the Liouville$-$Jacobi Identity\footnote{This result will be presented at the upcoming ``Third Workshop on Soliton Theory, Nonlinear Dynamics and Machine Learning" (August 18th$-$23rd, 2025) and the subsequent paper will be submitted for publication to the {\it Journal of Physics: Conference Series}}}
\author{Lubomir Markov 					\\
Department of Mathematics and CS 	\\
Barry University					\\
11300 N.E. Second Avenue		\\
Miami Shores, FL 33161, USA		\\
{\tt lmarkov@barry.edu}		\\}

\maketitle\noindent

\noindent
\text{MSC 2020: 34A30, 15A15} 
\bigskip

\begin{abstract} 
\noindent
A generalized Liouville$-$Jacobi Identity is proved for the determinant $\det{X(t)}$ of a solution $X(t)$ to the linear nonhomogeneous first-order matrix differential equation with left- and right-coefficient matrices $\ \frac{{\rm d}}{{\rm d}t} X(t) + A(t)X(t) + X(t)B(t)= F(t), \  X(t_0)=X_0.$

\end{abstract}
\bigskip\bigskip
\noindent
\section{Introduction}
Consider the linear nonhomogeneous first-order matrix differential equation involving left- and right-coefficient matrix multiplications:
\begin{equation}
\mathfrak{L}[X] \equiv \frac{{\rm d}}{{\rm d}t} X(t) + A(t)X(t) + X(t)B(t)= F(t), \quad X(t_0)=X_0.
\end{equation}

\smallskip

\noindent
Here $ \ X(t) \,=\,  \big[x_{ij}(t) \big]  \,=\, \big[{\bf x}_1 (t)  \ \, {\bf x}_2 (t)\ \,  \cdots \ \,  {\bf x}_n(t) \big]^{\rm T}, \quad A(t) = \big[a_{ij}(t) \big]$, \quad $B(t)= \big[b_{ij}(t)\big]\ $ and $\ F(t) = \big[f_{ij}(t) \big]  \,=\, \big[{\bf f}_1 (t)  \ \,   {\bf f}_2 (t)\ \,  \cdots \ \,  {\bf f}_n(t) \big]^{\rm T}$ are sufficiently smooth $n$-by-$n$ matrix-valued functions defined on some interval $[t_0, T\rangle,$ where $T$ can be finite or $+\infty$. The {\it homogeneous} case of this problem was introduced in a recent paper [3] wherein it was shown that if $X(t)$ satisfies $\mathfrak{L}[X] =0$, the following generalized Liouville's formula holds:

\begin{equation}
{\rm det} X(t) =  {\rm det} X_0 \,{\rm e}^{-\dint_{t_0}^t \{{\rm tr} \, A(t') + {\rm tr } \, B(t')\}\, {\rm d}t'} .
\end{equation}

When $B(t)=0$, formula (2) reduces to the classical identity due to Liouville (and Jacobi).
\noindent
In our previous work [5], we considered the nonhomogeneous differential equation 

\begin{equation}
\frac{{\rm d}}{{\rm d}t} X(t) + A(t)X(t) = F(t), \quad X(t_0)=X_0
\end{equation}
\noindent
and derived the following theorem which is another generalization of the Liouville$-$Jacobi Identity:

\begin{theorem}
The solution $X(t)$ to (3) satisfies
\begin{equation}
\quad {\rm det} X(t) \ = \ {\rm e}^{-\dint_{t_0}^t {\rm tr} \, A(\xi)\, {\rm d}\xi} \cdot \Bigg\{ \dint_{t_0}^t D_{X,F}(s)\, {\rm e}^{\, \int_{t_0}^s {\rm tr} \, A(\xi)\, {\rm d}\xi } \, {\rm d}s+  {\rm det} X_0 \Bigg\},
\end{equation}
where $ D_{X, F}(t) =\  \dsum_{j=1}^n {\rm det} \,\Big(X[{\bf x}_j \!\leftarrow {\bf f}_j ](t) \Big)$ and $\,X[{\bf x}_j \!\leftarrow {\bf f}_j ](t)\,$ is the matrix obtained from $X(t)$ by substituting the $j$th vector ${\bf x}_j (t)$ with the vector ${\bf f}_j (t) ,\  j=1, \dots, n. $ 
\end{theorem}

\bigskip\noindent
In what follows we shall derive an integral formula which generalizes all of the above-mentioned identities. Standard notation from linear algebra (see [1], [2]) and differential equations (see [6]) is used throughout, except for our denoting the adjugate of a matrix $X$ by $X^\#$ instead of the widespread $\rm{adj(X)}$ which we find very awkward. The following familiar properties (see [4]) of the adjugate matrix and the trace of a matrix shall be needed: 

\bigskip
$XX^\# = X^\# X = \det(X) I, \quad  \dfrac{{\rm d}}{{\rm d}t} \det X(t) = \rm{tr} \Big[X^\# \dfrac{{\rm d}X}{{\rm d}t}\Big] \text{\ (Jacobi's Formula)}, $
\medskip

$\rm{tr} (XY) = \rm{tr} (YX), \quad \rm{tr} (\lambda X + \mu Y) = \lambda\rm{tr}(X) + \mu\rm{tr}(Y) \text{\ (linearity of the trace).}$

\bigskip
\section{The General Liouville-Jacobi Identity}
\begin{theorem}
Suppose $X(t)$ solves the problem (1). Then for every $\,t\in [t_0, T\rangle$ we have

\begin{equation}
{\rm det} X(t)={\rm e}^{-\int_{t_0}^t {\rm tr}[A(\xi)+B(\xi)] {\rm d}\xi} \cdot \Big\{  \dint_{t_0}^t {\rm tr}[X^\#(s) F(s)]\, {\rm e}^{\, \int_{t_0}^s {\rm tr} [A(\xi)+B(\xi)] {\rm d}\xi } {\rm d}s+  {\rm det} X_0 \Big\}.
\end{equation} 
\end{theorem}

\noindent
{\bf Proof.}

The proof utilizes Jacobi's Formula in a powerful way. Beginning with it and performing standard manipulations, we have the following chain of equalities:
\begin{align}
	\dfrac{{\rm d}}{{\rm d}t} {\rm det} X(t) &=  {\rm tr} \Big[X^\#  \dfrac{{\rm d}X}{{\rm d}t}\Big]  =  
	{\rm tr} \big[X^\# (F -AX - XB) \big] 												\nonumber \\
		&= {\rm tr} \big[X^\# F \big] - {\rm tr} \big[X^\# AX \big] - {\rm tr} \big[X^\# X B\big]			\nonumber \\ 
		&= {\rm tr} \big[X^\# F \big] - {\rm tr} \big[AX X^\#\big] - {\rm tr} \big[X^\# X B\big]			\nonumber \\ 
		&= {\rm tr} \big[X^\# F \big] - {\rm tr} \big[A (\det {X}) I\big] - {\rm tr} \big[ (\det {X}) I B\big]	\nonumber \\ 
		&= {\rm tr} \big[X^\# F \big] - (\det {X}) {\rm tr} A  - (\det {X}) {\rm tr} B	= {\rm tr} \big[X^\# F\big] - 
				 {\rm tr}[A+B] \det {X}.												\nonumber \\ 
\nonumber 
\end{align}

\noindent
Thus the linear first-order scalar differential equation 
$$ \dfrac{{\rm d}}{{\rm d}t} {\rm det} X(t) +  {\rm tr}[A+B] \det {X} = {\rm tr} \big[X^\# F \big] $$ 
is obtained. Solving it by the standard application of an integrating factor yields (5) and concludes the proof.{\footnote{The possibility of using Jacobi's Formula is mentioned in [3] but not taken advantage of there.}} Note that (5) reduces to (4) when $B(t)=0,$ and it reduces to (2) when $F(t)=0.$

One more result is obtained with the additional assumption of the invertibility of the variables' matrix.
\begin{theorem}
Suppose that $X(t)$ solves the problem (1) and furthermore that it is invertible. Then for every $\,t\in [t_0, T\rangle$ we have
\begin{equation}
{\rm det} X(t)\ =\ {\rm det} X_0 \cdot {\rm e}^{\dint_{t_0}^t {\rm tr}\big[X^{-1}(s) F(s)-[A(s)+B(s)] \big] {\rm d}s} .
\end{equation} 
\end{theorem}

\noindent
{\bf Proof.}

If $X$ is invertible, Jacobi's Formula implies $\dfrac{{\rm d}}{{\rm d}t} \det X(t) = \det X(t) \,\rm{tr} \Big[X^{-1}\dfrac{{\rm d}X}{{\rm d}t}\Big],$
and (after manipulating $\rm{tr} \Big[X^{-1}\dfrac{{\rm d}X}{{\rm d}t}\Big]$ similarly to the previous proof) we obtain
\begin{equation}
\dfrac{\tfrac{{\rm d}}{{\rm d}t} \det X(t)}{ \det X(t)} = {\rm tr}[X^{-1}(t) F(t)] - {\rm tr}[A(t)+B(t)]. 
\end{equation} 
Integrating the last equation establishes the result.

\bigskip
All applications given in [3] of the identity (2) which corresponds to the homogeneous problem $\mathfrak{L}[X] =0$  will yield more general applications of the identity (5) corresponding to the nonhomogeneous problem $\mathfrak{L}[X] =F(t).$

\end{document}